\documentclass[]{article}
\usepackage{amssymb,amsfonts,amsmath}
\usepackage{aeguill}
\usepackage{textcomp}

\newtheorem{theo}{\sc Theorem}[section]

\newtheorem{lemma}[theo]{\sc Lemma}

\title{On special values of spinor $L$-functions of Siegel cusp eigenforms of genus $3$}
\author{Francesco Chiera 
\and Kirill Vankov}
\begin{document}
\maketitle

\abstract{We compute the special values for the spinor $L$-function $L(s,F_{12},spin)$ in the critical strip $s=12,\dotsc,19$, where $F_{12}$ is the unique (up to a scalar) Siegel cusp form of degree $3$ and weight $12$, which was constructed by Miyawaki.  These values are proportional to the product of Petersson inner products of Ramanujan's $\Delta$ by itself and the cusp form of weight 20 for $\mathrm{SL_2}(\mathbb{Z})$ by itself by a rational number and some power of $\pi$.  We also verify this result numerically using Dokchitser's ComputeL PARI package.  To our knowledge this is the first example of a spinor $L$-function of Siegel cusp forms of degree $3$, when the special values can be computed explicitly.}

\tableofcontents

\newpage
\section{Introduction}

In the remarkable paper \cite{miyawaki:1992} I.~Miyawaki considered certain Siegel cusp forms of degree 3, and on the basis of some numerical calculations, he was able to make interesting conjectures about the degeneration of the standard and spinor $L$-functions associated to such cusp forms.  Several years later T.~Ikeda \cite{ikeda:2006} proved Miyawaki's conjecture related to the standard $L$-function.  Basically he was able to construct an explicit lifting from Siegel cusp forms of degree $r$ to Siegel cusp forms of degree $r+2n$.  In particular, it turns out that the only cusp form of degree 3 and weight 12 is a basic example of this lifting (for $r=1$, $n=1$).

Recall that Miyawaki constructed his numerical examples by means of theta functions with spherical functions.  Namely, let $E_8$ be the unique even unimodular lattice of rank $8$ i.e.,
\begin{equation}
 E_8=\left\{ {}^t(x_1,\dotsc,x_8)\in\mathbb{R}^8 \left|\begin{aligned}&{2x_i\in\mathbb{Z} (i=1,\dotsc,8),}\\&{x_1+\dotsb+x_8\in2\mathbb{Z},}\\&{x_i-x_j\in\mathbb{Z}}\end{aligned}\right. \right\}\,,
\end{equation}
and
\begin{equation}
Q=\begin{pmatrix}1&0&0&i&0&0&0&0\\
                 0&1&0&0&i&0&0&0\\
                 0&0&1&0&0&i&0&0
       \end{pmatrix}
\end{equation}
be $3\times8$ matrix.  Then the theta series
\begin{equation}
\begin{split}
F_{12}&=\sum_{v_1,v_2,v_3\in E_8}\Re(\mathrm{det}\left(Q\cdot(v_1,v_2,v_3)\right)^8)\,\mathrm{exp}\left({\pi i\,\sigma((<v_i,v_j>)Z)}\right)\\
&=\sum_{N>0} a(N)\,\mathrm{exp}(2\pi i\,\sigma(NZ))
\end{split}
\end{equation}
is a cusp form of weight $12$ with respect to $\mathrm{Sp}_3(\mathbb{Z})$, where
\begin{equation}
\mathrm{Sp}_n(\mathbb{Z})=\left\{ M\in\mathrm{Mat}_{2n} \left| M\left(\begin{smallmatrix}0&\mathrm{I}_n\\-\mathrm{I}_n&0\end{smallmatrix}\right){}^tM=\left(\begin{smallmatrix}0&\mathrm{I}_n\\-\mathrm{I}_n&0\end{smallmatrix}\right) \right.\right\}\,.
\end{equation}

In the recent work \cite{heim:2008} B.~Heim proved Miyawaki's conjecture relevant to the spinor $L$-function for the specific case of the cusp form $F_{12}$.  In fact, Heim has showed that the following equality holds:
\begin{equation}
\label{eq:L(s,Sp(F_{12}))}
L(s,F_{12},spin)=L(s-9,\Delta)\,L(s-10,\Delta)\,L(s,\Delta\otimes g_{20})\,,
\end{equation}
where $\Delta$ is Ramanujan's discriminant cusp form and $g_{20}$ is the cusp form of weight $20$ of level $1$.  This theorem is the starting point of our investigation.  The purpose of this note is to show that one can provide the explicit rational number $R_s$ and power of $\pi$ such that for each critical value $s=12,\dotsc,19$
\begin{equation}
 L(s,F_{12},spin) = R_s\,\pi^{\alpha_s}\left<\Delta,\Delta\right>\,\left<g_{20},g_{20}\right>\,.
\end{equation}

We also compute the values numerically using the \textsc{SAGE} software \cite{SAGE} and Dokchitser's \textsc{ComputeL} \textsc{PARI} package \cite{ComputeL,dokchitser:2004}.

To our knowledge, this is the first example of a spinor $L$-function of Siegel cusp forms of degree $3$, when the special values can be computed explicitly.

It is possible to apply the same technique to compute the critical values of the spinor $L$-functions for non cuspidal modular forms; but there exist direct (and more simple) methods in this case.  For example, due to Zharkovskaya \cite[Theorem 1]{zarkovskaja:1974} there is the famous equality for non cuspidal forms
\begin{equation}
L(s,F,spin)=L(s,\Phi(F),spin)\,L(s-k+n,\Phi(F),spin)\,,
\end{equation}
where $\Phi$ is the Siegel operator.  We could also use our approach to compute the conjectural special values of the spinor $L$-function $L(s,F_{14},spin)$ for the unique Siegel cusp for of degree $3$ and weight $14$.

The authors are very grateful to Alexei Panchishkin for intensive and encouraging discussions, this work was undertaken due to his keen interest into the subject.  We also express our thanks to Bernhard Heim for useful remarks.

\section{Generalities and notations}

Let $\mathfrak{H}=\{z\in\mathbb{C}\mid\Im(z)>0\}$ be the upper half-plane.  For a positive integer $k$ and a Dirichlet character $\chi$ modulo a positive integer $N$ such that $\chi(-1)=(-1)^k$, we denote by $\mathcal{M}_k(\Gamma_0(N),\chi)$ the vector space of all holomorphic modular forms $f(z)$ of weight $k$ satisfying
\begin{equation}
f(\gamma(z))=\chi(d)(cz+d)^kf(z)\text{~for~all~}\gamma=\begin{pmatrix}a&b\\c&d\end{pmatrix}\in\Gamma_0(N)\,,
\end{equation}
where the variable $z\in\mathfrak{H}$, $\displaystyle\gamma(z)=\frac{az+b}{cz+d}$, and
\begin{equation}
\Gamma_0(N)=\left\{\begin{pmatrix}a&b\\c&d\end{pmatrix}\in\mathrm{SL}_2(\mathbb{Z})\mid c\equiv 0(\mathrm{mod}\,N)\right\}\,.
\end{equation}
We denote by $S_k(N,\chi)$ the subspace of $\mathcal{M}_k(\Gamma_0(N),\chi)$ consisting of all cusp forms.  Every element $f$ of $\mathcal{M}_k(\Gamma_0(N),\chi)$ has a Fourier expansion
\begin{equation}
f(z)=\sum_{n=0}^\infty a(n)\,q^n\,,
\end{equation}
where $q=\exp(2\pi iz)$ and $a(n)$ are complex numbers in general.

The $L$-function associated to $f$ is defined as $\displaystyle L(s,f)=\sum_{n=1}^\infty a(n)\,n^{-s}$.  More generally, with an arbitrary Dirichlet character $\omega$, the twisted $L$-function is defined as $\displaystyle L(s,f,\omega)=\sum_{n=1}^\infty a(n)\,\omega(n)\,n^{-s}$.  These $L$-functions can be also written in the form of Euler product:
\begin{align}
L(s,f)&=\prod_{p\text{ prime}}(1-a(p)\,p^{-s}+\chi(p)\,p^{k-1-2s})^{-1}\,,\\
L(s,f,\omega)&=\prod_{p\text{ prime}}(1-a(p)\,\omega(p)\,p^{-s}+\chi(p)\,\omega(p)^2\,p^{k-1-2s})^{-1}\,.
\end{align}
Next, the Dirichlet $L$-series for any character $\chi$ of conductor $N$ is given by 
\begin{equation}
\label{eq:Dirichlet_L-series}
L(s,\chi)=\sum_{n=1}^\infty\chi(n)\,n^{-s}
\end{equation}
and its Euler product 
\begin{equation}
\label{eq:Euler_product_Dirichlet_L-series}
L(s,\chi)=\prod_{p\,\nmid N}\frac{1}{1-\chi(p)\,p^{-s}}\,.
\end{equation}
In the case, when $\chi$ is the identity Dirichlet character, the latter series is Riemann's zeta function $\displaystyle\zeta(s)=\sum_{n=1}^\infty\frac{1}{n^s}$.

Let $\displaystyle g(z)=\sum_{n=0}^\infty b(n)\,q^n \in \mathcal{M}_l(\Gamma_0(N),\xi)$ be another modular form of weight $l$ with Fourier coefficients $b(n)$.  The $L$ function associated to two modular forms $f$ and $g$ is given by the additive convolution
\begin{equation}
L(s,f,g)=\sum_{n=1}^\infty a(n)\,b(n)\,n^{-s}\,.
\end{equation}

Another type of $L$-functions associated to two modular forms is Rankin's product $L$-function (multiplicative convolution).  It is denoted by \mbox{$L(s,f\otimes g)$} and defined as (see \cite[page 786]{shimura:1976}):
\begin{equation}
\label{eq:Rankins_identity}
L(s,f\otimes g)=L_N(2s+2-k-l,\chi\xi)\,L(s,f,g)\,,
\end{equation}
where $L_N(s,\omega)$ with a Dirichlet character $\omega$ modulo $N$ is defined, as usual, by (\ref{eq:Dirichlet_L-series}) with $\omega(n)=0$ for $(n,N)\neq 1$, and the Euler factors in (\ref{eq:Euler_product_Dirichlet_L-series}), corresponding to the prime divisors of a number $N$, have been omitted.  Note, that in the case, when $f$ and $g$ are cusp eigenforms, the lefthand side of (\ref{eq:Rankins_identity}) is an Euler product of degree $4$ in view of the following Lemma:
\begin{lemma}[Lemma 1, \cite{shimura:1976}]
\label{lemma:Rankin}
Suppose we have formally
$$
\sum_{n=1}^\infty A(n)\,n^{-s}=\prod_p\left[(1-\alpha_p\,p^{-s})(1-\alpha_p^\prime\,p^{-s})\right]^{-1}
$$
and
$$
\sum_{n=1}^\infty B(n)\,n^{-s}=\prod_p\left[(1-\beta_p\,p^{-s})(1-\beta_p^\prime\,p^{-s})\right]^{-1}\,.
$$
Then 
$$
\sum_{n=1}^\infty A(n)\,B(n)\,n^{-s}=\frac{1-\alpha\alpha^\prime\beta\beta^\prime p^{-2s}}{(1-\alpha\beta p^{-s})(1-\alpha\beta^\prime p^{-s})(1-\alpha^\prime\beta p^{-s})(1-\alpha^\prime\beta^\prime p^{-s})}\,.
$$
\end{lemma}

For two elements $f,h\in\mathcal{M}_k(\Gamma_0(N))$ such that $fh$ is a cusp form, the Petersson inner product $\left<f,h\right>$ is defined as
\begin{equation}
\label{eq:Petersson}
\left<f,h\right>=\dfrac{1}{[\mathrm{SL}_2(\mathbb{Z}):\Gamma_0(N)]}\int_{\Phi_N}\overline{f(z)}\,h(z)\,y^{k-2}\,dx\,dy\,,
\end{equation}
where $z=x+iy$, $\Phi_N$ is a fundamental domain for $\mathfrak{H}$ modulo $\Gamma_0(N)$ and the bar denotes the complex conjugate.  We also define $\left<f,h\right>$ by (\ref{eq:Petersson}) for nearly holomorphic modular forms $f$ and $h$ on $\mathfrak{H}$ whenever the integral is convergent (see \cite[section 8.2]{shimura:2007} for definition and properties of the nearly holomorphic modular forms).

We also give a brief definition of the spinor $L$-function.  For a Hecke eigenform $F$ with local Satake parameters $\mu_0,\mu_1,\dotsc,\mu_n$ it is given by the infinite product 
\begin{equation}
 L(s,F,spin)=\prod_pL_p(s,F,spin),
\end{equation}
where for each prime number $p$
\begin{equation}
 L_p(s,F,spin)=\left((1-\mu_0\,p^{-s})\prod_{r=1}^n\prod_{i_1<\dotsb<i_r}(1-\mu_0\,\mu_{i_1}\,\dotsb\mu_{i_r}\,p^{-s})\right)^{-1}\,.
\end{equation}
We do not use this definition in our computation, but develop the righthand side of the identity (\ref{eq:L(s,Sp(F_{12}))}).  The explicit description of the Satake parameters for $F_{12}$ is given in \cite{heim:2008}.

\section{The expression for $L(s,\Delta)\,L(s-1,\Delta)$}

Let $f=\Delta$ be Ramanujan's discriminant modular form of weight $k=12$:
\begin{equation}
\begin{split}
f(z)&=\Delta(z)=\sum_{n=1}^\infty\tau(n)\,q^n=q\prod_{n=1}^\infty\left(1-q^n\right)^{24}\\
&=q-24\,q^2+252\,q^3-1472\,q^4+4830\,q^5-6048\,q^6+\cdots\,,
\end{split}
\end{equation}
where $\tau(n)$ is Ramanujan's tau function.  The associated $L$-function is 
\begin{equation}
L(s,f)=\sum_{n=1}^\infty\tau(n)\,n^{-s}=\prod_{p}\left(1-\tau(p)\,p^{-s}+p^{11-2s}\right)^{-1}\,.
\end{equation}
Let $\displaystyle G_2(z)=-\frac{1}{24}+\sum_{n=1}^\infty\sigma_1(n)\,q^n=-\frac{1}{24}+q+3\,q^2+4\,q^3+7\,q^4+\cdots$, where $\displaystyle\sigma_1(n)=\sum_{d|n}d$ ~is the divisor function defined as the sum of the divisors of~$n$.  Consider the Eisenstein series (see \cite[Lemma 7.2.19, (2)]{miyake:2006})
\begin{equation}
g(z)=G_{2,p}(z)=G_2(z)-p\,G_2(pz)=\frac{p-1}{24}+\sum_{n=1}^\infty\sum_{\substack{d|n\\p\,\nmid\,d}}d\,q^n\,,
\end{equation}
of weight $2$ for $\Gamma_0(p)$ and the corresponding $L$ series
\begin{equation}
\begin{split}
L(s,g)&=\sum_{n=1}^\infty\sum_{\substack{d|n\\p\,\nmid\,d}}d\,n^{-s}=\sum_{\substack{d,d_1\geqslant 1\\p\,\nmid\,d}}d\,(d\,d_1)^{-s}=\sum_{\substack{d\geqslant 1\\p\,\nmid\,d}}d^{1-s}\sum_{d_1\geqslant 1}d_1^{-s}\\
\qquad&=(1-p^{1-s})\,\zeta(s-1)\,\zeta(s)\,.
\end{split}
\end{equation}
Let us put $p=2$, then 
\begin{equation}
\label{eq:G_2_2}
\begin{split}
g(z)&=G_{2,2}(z)=G_2(z)-2\,G_2(2z)\\
&=\frac{1}{24}+q+q^2+4\,q^3+q^4+6\,q^5+4\,q^6+\cdots\,,
\end{split}
\end{equation}
and
\begin{equation}
L(s,g)=(1-2^{1-s})\,\zeta(s-1)\,\zeta(s)\,.
\end{equation}
For $f=\Delta=\sum\tau(n)\,q^n\in S_{12}(2)$ and $g=G_{2,2}=\sum b(n)\,q^n\in\mathcal{M}_2(\Gamma_0(2),\xi)$ we have $k=12$ and $l=2$.  Put $N=2$, $\chi=1$ and 
\begin{equation}
\xi(n)=(1~\mathrm{mod}~N)(n)=\begin{cases}1,&\text{if $n$ odd;}\\0,&\text{if $n$ even.}\end{cases}
\end{equation}
Assume $1-\tau(p)\,X+p^{11}\,X^2=(1-\alpha_p\,X)(1-\alpha_p^\prime\,X)$, then $\alpha_p+\alpha_p^\prime=\tau(p)$, $\alpha_p\,\alpha_p^\prime=p^{11}$, and 
\begin{equation}
L(s,f)=\prod_{p}\left((1-\alpha_p\,p^{-s})(1-\alpha_p^\prime\,p^{-s})\right)^{-1}\,.
\end{equation}
Similarly, consider 
\begin{equation}
\begin{split}
L(s,g)&=\sum_{n=1}^\infty b(n)\,n^{-s}\\
&=(1-2^{1-s})\,\zeta(s-1)\,\zeta(s)\\
&=(1-2^{1-s})\prod_p\left((1-p^{1-s})\,(1-p^{-s})\right)^{-1}\\
&=\prod_p\left((1-\beta_p\,p^{-s})\,(1-\beta_p^\prime\,p^{-s})\right)^{-1}\,,
\end{split}
\end{equation}
where $\beta(p)=1$ for all $p$, $\beta^\prime(2)=0$ and $\beta^\prime(p)=p$ for all odd primes.  By definition and using Lemma \ref{lemma:Rankin} 
\begin{equation}
\label{eq:L_Rankin}
\begin{split}
&L(s,\Delta\otimes G_{2,2})=L_2(2s+2-12-2,\psi)\,L(s,f,g)\\
&=\prod_{p\neq 2}(1-p^{12-2s})^{-1}\,\cdot\,\sum_{n=1}^\infty\tau(n)\,b(n)\,n^{-s}\\
&=\prod_{p\neq 2}(1-p^{12-2s})^{-1}\times\\
&\quad\times\prod_p\frac{1-\alpha_p\,\alpha_p^\prime\,\beta_p\,\beta_p^\prime\,p^{-2s}}{(1-\alpha_p\,\beta_p\,p^{-s})\,(1-\alpha_p^\prime\,\beta_p\,p^{-s})\,(1-\alpha_p\,\beta_p^\prime\,p^{-s})\,(1-\alpha_p^\prime\,\beta_p^\prime\,p^{-s})}\\
&=\frac{1}{(1-\alpha_2\,2^{-s})\,(1-\alpha_2^\prime\,2^{-s})}\times\\
&\quad\times\prod_{p\neq2}\frac{(1-p^{12-2s})^{-1}\,(1-p^{11}\,p\,p^{-2s})}{(1-\alpha_p\,p^{-s})\,(1-\alpha_p^\prime\,p^{-s})\,(1-\alpha_p\,p\,p^{-s})\,(1-\alpha_p^\prime\,p\,p^{-s})}\\
&=\frac{1}{(1-\alpha_2\,2^{-s})\,(1-\alpha_2^\prime\,2^{-s})}\times\\
&\quad\times\prod_{p\neq2}\frac{1}{\left((1-\alpha_p\,p^{-s})\,(1-\alpha_p^\prime\,p^{-s})\right)\,\left((1-\alpha_p\,p^{1-s})\,(1-\alpha_p^\prime\,p^{1-s})\right)}\\
&=\prod_p\frac{1}{(1-\alpha_p\,p^{-s})\,(1-\alpha_p^\prime\,p^{-s})}\,\prod_{p\neq2}\frac{1}{(1-\alpha_p\,p^{1-s})\,(1-\alpha_p^\prime\,p^{1-s})}\\
&=L(s,\Delta)\,L(s-1,\Delta)\,(1-\tau(2)\,2^{1-s}+2^{11}\,s^{2-2s})\\
&=(1+3\cdot2^{4-s}+s^{13-2s})\,L(s,\Delta)\,L(s-1,\Delta)\,.
\end{split}
\end{equation}

Finally, we obtain the following identity
\begin{equation}
\label{eq:L(s,Delta)L(s-1,Delta)}
L(s,\Delta)\,L(s-1,\Delta)=\frac{L(s,\Delta\otimes G_{2,2})}{(1+3\cdot2^{4-s}+2^{13-2s})}\,.
\end{equation}

\section{Computation of $L(s,\Delta\otimes G_{2,2})$}
\label{sec:L(s,Delta_x_G(2,2))}

Now we express $L(s,\Delta\otimes G_{2,2})$ (at integral points $s=3,\dotsc,10$) as a multiple of Petersson inner product $\left<\Delta,\Delta\right>$.  Using \cite[(2.4)]{shimura:1976},
\begin{equation}
\label{eq:Shimura-2.4}
\begin{split}
&L(s,\Delta\otimes G_{2,2})\\
&=\frac{(4\,\pi)^s}{2\,\Gamma(s)}\int_{\Phi_2}\overline{\Delta(z)}\,G_{2,2}(z)\,E_{10,2}(z,s-11,\xi)\,y^{s-1}\,dx\,dy\\
&=\frac{(4\,\pi)^s}{2\,\Gamma(s)}\int_{\Phi_2}\overline{\Delta(z)}\,G_{2,2}(z)\,E_{10,2}(z,s-11,\xi)\,y^{s-11}\,y^{10}\,dx\,dy\\
&=\frac{(4\,\pi)^s\,[\mathrm{SL}_2(\mathbb{Z}):\Gamma_0(2)]}{2\,\Gamma(s)}\left<\Delta(z),G_{2,2}(z)\,y^{s-11}\,E_{10,2}(z,s-11,\xi)\right>\\
&=\frac{3}{2}\frac{(4\,\pi)^s}{\Gamma(s)}\left<\Delta(z),\mathcal{H}ol\left(G_{2,2}(z)\,y^{s-11}\,E_{10,2}(z,s-11,\xi)\right)\right>\\
&=\frac{3}{2}\frac{(4\,\pi)^{11}}{\Gamma(s)}\left<\Delta(z),\mathcal{H}ol\left(G_{2,2}(z)\,(4\pi y)^{s-11}\,E_{10,2}(z,s-11,\xi)\right)\right>\,,
\end{split}
\end{equation}
where $\left<f,g\right>$ is the Petersson inner product (\ref{eq:Petersson}), $\Phi_2$ denotes a fundamental domain for $\Gamma_0(2)\backslash\mathfrak{H}$, $z=x+iy$, 
\begin{equation}
E_{\lambda,N}(z,s,\xi)=\mathop{{\sum}^\prime}_{(m,n)} \xi(n)\,(mNz+n)^{-\lambda}|mNz+n|^{-2s}\,,
\end{equation}
$\mathop{{\sum}^\prime}$ denotes the summation over all $(m,n)\in\mathbb{Z}^2$, $(m,n)\neq(0,0)$, $\mathcal{H}ol(F)$ is the operator of holomorphic projection  (see \cite[(2.148)]{courtieu-panchishkin:2004}).  It is defined so, that $\left<f,F\right>=\left<f,\mathcal{H}ol(F)\right>$ for all $f\in S_k(N,\psi)$.

In order to compute the holomorphic projection of the product in the last identity of (\ref{eq:Shimura-2.4}) $\mathcal{H}ol\left(G_{2,2}(z)\,(4\pi y)^{s-11}\,E_{10,2}(z,s-11,\xi)\right)$ we write the Fourier expansion for $E_{10,2}(z,s-11,\xi)$ in the convenient form using the Whittaker functions by applying the Proposition 2.2 of \cite{panchishkin:2003}, where in notations of that proposition 
\begin{equation*}
\begin{split}
&s\in\{3,\dotsc,10\}\,,\quad s-11<0\,,\quad s-11+10=s-1>0\,,\\
&a=0\,,\quad b=1\,,\\
&\mathbf{E}_{10,2}(z,s-11,0,1)=E_{10,2}(z,s-11,\xi)\,,\\
&\delta(\frac{a}{N})=1\,,\\
&\zeta(s;0,2)=\sum_{0<n\equiv0(2)}n^{-s}=\sum_{n=1}^\infty(2n)^{-s}=2^{-s}\zeta(s)\,,\\
&\zeta(s;1,2)=\sum_{0<n\equiv1(2)}n^{-s}=\zeta(s)-\zeta(s;0,2)=(1-2^{-s})\zeta(s)\,.
\end{split}
\end{equation*}
The Whittaker function $W(y,\alpha,\beta)$ is defined as (see \cite[(2.4)]{panchishkin:2003}, for example)
\begin{equation}
W(y,\alpha,\beta)=\Gamma(\beta)^{-1}\int_0^{+\infty}(u+1)^{\alpha-1}\,u^{\beta-1}\,\mathrm{e}^{-yu}\,du
\end{equation}
for $y>0$, $\alpha,\beta\in\mathbb{C}$ with $\Re(\beta)>0$ and for arbitrary $\alpha$ and $\beta$ this function is defined by the analytic continuation and the functional equation:
\begin{equation}
W(y,\alpha,\beta)=y^{1-\alpha-\beta}\,W(y,1-\beta,1-\alpha)\,.
\end{equation}
For a non negative integer $r$, we have 
\begin{equation}
W(y,\alpha,-r)=\sum_{i=0}^{r}\frac{(-1)^i\binom{r}{i}\Gamma(\alpha)}{\Gamma(\alpha-i)}\,y^{r-i}\,.
\end{equation}

Therefore,
\begin{equation}
\begin{split}
&(4\pi y)^{s-11}\,E_{10,2}(z,s-11,\xi)\\
&=(4\pi y)^{s-11}\,2\,\zeta(2s-12;1,2)\\
&\quad+(4\pi y)^{s-11}\,\frac{(-2\pi i)^{2s-12}\,(-1)^{s-11}\Gamma(2s-13)}{(4\pi y)^{2s-13}\,2\,\Gamma(s-1)\,\Gamma(s-11)}\,2\,\zeta(2s-13;0,2)\\
&\quad+(4\pi y)^{s-11}\,\frac{(-2\pi i)^{2s-12}\,(-1)^{s-11}}{2^{2s-12}\Gamma(s-1)}\times\\
&\qquad\times \sum_{n=1}^\infty\sum_{\pm d|n}\mathrm{sgn}(d)\,d^{\,2s-13}\,\mathrm{e}^{\pi id}\,W(4\pi ny,s-1,s-11)\,q^n\\
&=2\,(4\pi y)^{s-11}\,(1-2^{12-2s})\,\zeta(2s-12)\\
&\quad-(4\pi y)^{2-s}\frac{2\pi^{2s-12}\,\Gamma(2s-13)\,\zeta(2s-13)}{\Gamma(s-11)\,\Gamma(s-1)}\\
&\quad-(4\pi y)^{s-11}\frac{2\pi^{2s-12}}{\Gamma(s-1)}\sum_{n=1}^\infty\sum_{d|n}(-1)^dd^{\,2s-13}\,W(4\pi ny,s-1,s-11)\,q^n\,.
\end{split}
\end{equation}
Let 
\begin{equation}
\begin{split}
C_0^\prime&=C_0^\prime(s)=(-1)\frac{2\,\pi^{2s-12}\,\Gamma(2s-13)\,\zeta(2s-13)}{\Gamma(s-11)\,\Gamma(s-1)}\,,\\
C_0^{\prime\prime}&=C_0^{\prime\prime}(s)=(2-2^{13-2s})\,\zeta(2s-12)\,,\\
C_1&=C_1(s)=2\,\pi^{2s-12}\,,\\
C_2&=C_2(s)=(2-2^{2s-12})\,\pi^{2s-12}\,,
\end{split}
\end{equation}
then 
\begin{equation}
\begin{split}
(4\pi y)^{s-11}&\,E_{10,2}(z,s-11,\xi)\\
&=C_0^\prime\,(4\pi y)^{2-s}+C_0^{\prime\prime}\,(4\pi y)^{s-11}\\
&\quad+C_1\,\frac{W(4\pi y,s-1,s-11)}{\Gamma(s-1)}(4\pi y)^{s-11}\,q\\
&\quad+C_2\,\frac{W(8\pi y,s-1,s-11)}{\Gamma(s-1)}(4\pi y)^{s-11}\,q^2+\cdots\,.
\end{split}
\end{equation}
Recall that $G_{2,2}(z)=1/24+q+q^2+\cdots$ by (\ref{eq:G_2_2}).  We write the Fourier coefficients $\widetilde{A}_i(s,y)$ for the product $F=G_{2,2}(z)\,(4\pi y)^{s-11}\,E_{10,2}(z,s-11,\xi)=\sum \widetilde{A}_n(s,y)q^n$ in order to apply the Holomorphic Projection Lemma \cite[Proposition (5.1)]{gross-zagier:1986} to find the image of the projection operator $\mathcal{H}ol(F)=\sum A_n(s)q^n$ (the Lemma is originally due to Sturm \cite{sturm:1980}).  It should be noted, that the relevant polynomial decay hypotheses of the Lemma are satisfied for all actions $E_{10,2}(z,s-11,\xi)\vert\gamma$ of $\gamma\in\mathrm{SL}_2(\mathbb{Z})$ and each critical point $s$, see \cite[(2.3)]{panchishkin:2003}.

We need to compute just two coefficients $A_1=A_1(s)$ and $A_2=A_2(s)$ since the result of the holomorphic projection belongs to the space of cusp forms for the subgroup $\Gamma_0(2)$, which has the dimension $2$.  Then we find the linear combination representing $\mathcal{H}ol(F)$ in the basis $\{\Delta(z),\Delta(2z)\}$:
\begin{equation}
\label{eq:Hol=alpha*beta}
\mathcal{H}ol(F)=\alpha\cdot\Delta(z)+\beta\cdot\Delta(2z)\,.
\end{equation}
Namely, the computation of $A_1$ and $A_2$ gives:
\begin{equation}
\label{eq:A1}
\begin{split}
&A_1 = \frac{C_0^\prime}{10!}\int_0^{+\infty}(4\pi y)^{2-s}\,\mathrm{e}^{-4\pi y}\,(4\pi y)^{10}\,d(4\pi y)\\
&\quad +\frac{C_0^{\prime\prime}}{10!}\int_0^{+\infty}(4\pi y)^{s-11}\,\mathrm{e}^{-4\pi y}\,(4\pi y)^{10}\,d(4\pi y)\\
&\quad +\frac{C_1}{24\cdot 10!}\int_0^{+\infty}\frac{W(4\pi y,s-1,s-11)}{\Gamma(s-1)}(4\pi y)^{s-11}\,\mathrm{e}^{-4\pi y}\,(4\pi y)^{10}\,d(4\pi y)\\
&=\frac{\Gamma(13-s)}{10!}\,C_0^\prime + \frac{\Gamma(s)}{10!}\,C_0^{\prime\prime}\\
&\quad +\frac{C_1}{24\cdot 10!}\int_0^{+\infty}\sum_{i=0}^{11-s}\frac{(-1)^i\binom{11-s}{i}(4\pi y)^{10-i}}{\Gamma(s-1-i)}\,\mathrm{e}^{-4\pi y}\,d(4\pi y)\\
&=\frac{\Gamma(13-s)}{10!}\,C_0^\prime + \frac{\Gamma(s)}{10!}\,C_0^{\prime\prime}\\
&\quad +\frac{C_1}{24\cdot 10!}\sum_{i=0}^{11-s}\frac{(-1)^i\binom{11-s}{i}\Gamma(11-i)}{\Gamma(s-1-i)}\,,
\end{split}
\end{equation}
\begin{equation}
\label{eq:A2}
\begin{split}
&A_2 = \frac{C_0^\prime}{10!}\int_0^{+\infty}(4\pi y)^{2-s}\,\mathrm{e}^{-8\pi y}\,(8\pi y)^{10}\,d(8\pi y)\\
&\quad +\frac{C_0^{\prime\prime}}{10!}\int_0^{+\infty}(4\pi y)^{s-11}\,\mathrm{e}^{-8\pi y}\,(8\pi y)^{10}\,d(8\pi y)\\
&\quad +\frac{C_1}{10!}\int_0^{+\infty}\frac{W(4\pi y,s-1,s-11)}{\Gamma(s-1)}(4\pi y)^{s-11}\,\mathrm{e}^{-8\pi y}\,(8\pi y)^{10}\,d(8\pi y)\\
&\quad +\frac{C_2}{24\cdot 10!}\int_0^{+\infty}\frac{W(8\pi y,s-1,s-11)}{\Gamma(s-1)}(4\pi y)^{s-11}\,\mathrm{e}^{-8\pi y}\,(8\pi y)^{10}\,d(8\pi y)\\
&=\frac{C_0^\prime}{10!}\,2^{s-2}\int_0^{+\infty}(8\pi y)^{12-s}\,\mathrm{e}^{-8\pi y}\,d(8\pi y)\\
&\quad +\frac{C_0^{\prime\prime}}{10!}\,2^{11-s}\int_0^{+\infty}(8\pi y)^{s-1}\,\mathrm{e}^{-8\pi y}\,d(8\pi y)\\
&\quad +\frac{C_1}{10!}\int_0^{+\infty}\sum_{i=0}^{11-s}\frac{(-1)^i\binom{11-s}{i}(4\pi y)^{11-s-i}}{\Gamma(s-1-i)}\,(4\pi y)^{s-11}\,\mathrm{e}^{-8\pi y}\,(8\pi y)^{10}\,d(8\pi y)\\
&\quad +\frac{C_2}{24\cdot 10!}\int_0^{+\infty}\sum_{i=0}^{11-s}\frac{(-1)^i\binom{11-s}{i}(8\pi y)^{11-s-i}}{\Gamma(s-1-i)}\,(4\pi y)^{s-11}\,\mathrm{e}^{-8\pi y}\,(8\pi y)^{10}\,d(8\pi y)\\
&=\frac{\Gamma(13-s)}{10!}\,2^{s-2}\,C_0^\prime + \frac{\Gamma(s)}{10!}\,2^{11-s}\,C_0^{\prime\prime}\\
&\quad +\frac{C_1}{10!}\sum_{i=0}^{11-s}2^i\,\frac{(-1)^i\binom{11-s}{i}\Gamma(11-i)}{\Gamma(s-1-i)}\\
&\quad +\frac{C_2}{24\cdot 10!}\,2^{11-s}\,\sum_{i=0}^{11-s}\frac{(-1)^i\binom{11-s}{i}\Gamma(11-i)}{\Gamma(s-1-i)}\,.
\end{split}
\end{equation}

These Fourier coefficients are rational numbers up to the factor $\pi^{2s-12}$ for each $s\in\{3,\dotsc,10\}$, see Table \ref{tbl:A1A2}.

\begin{table}[!tb]
\caption{Fourier coefficient of $\mathcal{H}ol\left(G_{2,2}(z)\,(4\pi y)^{s-11}\,E_{10,2}(z,s-11,\xi)\right)$}\label{tbl:A1A2}
$$
\renewcommand{\arraystretch}{0}
\begin{array}{|c|c|c|c|}\hline
\rule{0pt}{11pt}s&\pi\text{-factor}&A_1(s)&A_2(s)\\
\rule{0pt}{2pt}&&&\\\hline\rule{0pt}{2pt}&&&\\
\strut  3&\pi^{-6}&\dfrac{1}{50}&\dfrac{76}{25}\\
\rule{0pt}{2pt}&&&\\\hline\rule{0pt}{2pt}&&&\\
\strut  4&\pi^{-4}&\dfrac{1}{270}&\dfrac{56}{135}\\
\rule{0pt}{2pt}&&&\\\hline\rule{0pt}{2pt}&&&\\
\strut  5&\pi^{-2}&\dfrac{1}{1440}&\dfrac{1}{20}\\
\rule{0pt}{2pt}&&&\\\hline\rule{0pt}{2pt}&&&\\
\strut  6&\pi^{ 0}&\dfrac{1}{6048}&\dfrac{5}{756}\\
\rule{0pt}{2pt}&&&\\\hline\rule{0pt}{2pt}&&&\\
\strut  7&\pi^{ 2}&\dfrac{1}{16800}&\dfrac{1}{900}\\
\rule{0pt}{2pt}&&&\\\hline\rule{0pt}{2pt}&&&\\
\strut  8&\pi^{ 4}&\dfrac{17}{518400}&\dfrac{13}{64800}\\
\rule{0pt}{2pt}&&&\\\hline\rule{0pt}{2pt}&&&\\
\strut  9&\pi^{ 6}&\dfrac{11}{453600}&\dfrac{-1}{56700}\\
\rule{0pt}{2pt}&&&\\\hline\rule{0pt}{2pt}&&&\\
\strut 10&\pi^{ 8}&\dfrac{13}{604800}&\dfrac{-1}{10800}\\
\rule{0pt}{2pt}&&&\\\hline
\end{array}
$$
\end{table}

Next, we compute $\alpha$ and $\beta$ in the linear combination (\ref{eq:Hol=alpha*beta}) by comparing the coefficients $A_1$ and $A_2$ corresponding to terms $q$ and $q^2$ with the equivalent linear combination of coefficients of our basis functions $\Delta(z)=q-24q^2+\dotsc$ and $\Delta(2z)=q^2+\dotsc$:
\begin{equation}
\left\{
\begin{aligned}
A_1 &= \alpha\cdot1 + \beta\cdot0\\
A_2 &= \alpha\cdot(-24) + \beta\cdot1
\end{aligned}
\right.
\end{equation}
Resolving this system of linear equations, we obtain 
\begin{equation}
\begin{split}
\alpha &= A_1\\
\beta &= 24\,A_1+A_2\,,
\end{split}
\end{equation}
therefore, we obtain the following identity for the Rankin's convolution of $\Delta$ and $G_{2,2}$ (\ref{eq:Shimura-2.4}):
\begin{equation}
\label{eq:delta(z)*delta(2z)}
\begin{split}
L(s&,\Delta\otimes G_{2,2})\\
&=\frac{3}{2}\frac{(4\pi)^{11}}{\Gamma(s)}\left(A_1(s)\left<\Delta(z),\Delta(z)\right>+(24\,A_1(s)+A_2(s))\left<\Delta(z),\Delta(2z)\right>\right).
\end{split}
\end{equation}
% Since the group index $[\mathrm{SL}_2(\mathbb{Z}):\Gamma_0(2)]$ is equal to $3$, we can write
% \begin{equation}
% \label{eq:delta(z)*delta(2z)}
% \begin{split}
% L(s&,\Delta\otimes G_{2,2})\\
% &=\frac{(4\pi)^{11}}{6\Gamma(s)}\left(A_1\left<\Delta(z),\Delta(z)\right>+(24\,A_1+A_2)\left<\Delta(z),\Delta(2z)\right>\right)\,,
% \end{split}
% \end{equation}
% where we drop the subscript in Petersson inner product $\left<\cdot,\cdot\right>$, since now it is evaluated for the full modular group $\mathrm{SL}_2(\mathbb{Z})$.
We simplify the obtained expression even further.  Recall that 
\begin{equation}
 \Delta(2z)=2^{-k/2}\,\Delta(z)\vert V(2)\,,
\end{equation}
where $k=12$ is the weight of $\Delta$ and $V(2)$ is the operator $\begin{pmatrix}2&0\\0&1\end{pmatrix}$.  Then 
\begin{equation}
 \left<\Delta(z),\Delta(2z)\right>=2^{-6}\left<\Delta(z),\Delta(z)\vert\left(\begin{smallmatrix}2&0\\0&1\end{smallmatrix}\right)\right>\,.
\end{equation}
Consider $\gamma\in\Gamma_0(2)\backslash\Gamma$, $\Gamma=\mathrm{SL}_2(\mathbb{Z})$.  The summation over all $\gamma$ gives 
\begin{equation}
\begin{split}
 \left<\Delta(z),\Delta(2z)\right>&=2^{-6}[\Gamma:\Gamma_0(2)]^{-1}\sum_\gamma\left<\Delta(z)\vert\gamma,\Delta(z)\vert\left(\begin{smallmatrix}2&0\\0&1\end{smallmatrix}\right)\gamma\right>\\
&=2^{-6}3^{-1}\left<\Delta(z),\mathrm{Tr^{(2)}}\left(\Delta(z)\vert\left(\begin{smallmatrix}2&0\\0&1\end{smallmatrix}\right)\right)\right>\,.
\end{split}
\end{equation}
The trace operator $\mathrm{Tr^{(N)}}:\mathcal{M}_k(\Gamma_0(N))\rightarrow\mathcal{M}_k(\mathrm{SL}_2(\mathbb{Z}))$ is defined by the action $f\rightarrow\displaystyle\sum_{\gamma\in\Gamma_0(N)\backslash\mathrm{SL}_2(\mathbb{Z})}f|_k\gamma$.  We have (see \cite{serre:1973})
\begin{equation}
\mathrm{Tr^{(2)}}\left(\Delta(z)\vert\left(\begin{smallmatrix}2&0\\0&1\end{smallmatrix}\right)\right)=2^{-5}\,T_2(\Delta)\,
\end{equation}
% The action of the Hecke operator $T_2$ can be written as 
% \begin{equation}
%  T_2(f) = 2^{-\frac{k}{2}+1}\sum_{\delta\in\Gamma\backslash\Gamma\left(\begin{smallmatrix}2&0\\0&1\end{smallmatrix}\right)\Gamma}f\vert\delta\,,
% \end{equation}
where $T_2$ is the Hecke operator, therefore, 
\begin{equation}
\begin{split}
 \left<\Delta(z),\Delta(2z)\right>&=2^{-6}3^{-1}\left<\Delta(z),2^{-5}\,T_2\left(\Delta(z)\right)\right>\\
 &=2^{-11}3^{-1}\left<\Delta(z),\tau(2)\,\Delta(z)\right>\\
 &=-\frac{1}{256}\left<\Delta(z),\Delta(z)\right>\,.
\end{split}
\end{equation}
Substituting the last identity into (\ref{eq:delta(z)*delta(2z)}), we obtain the final expression:
\begin{equation}
\label{eq:part1}
L(s,\Delta\otimes G_{2,2})=\frac{3}{2}\frac{(4\pi)^{11}}{\Gamma(s)}\frac{(232\,A_1(s)-A_2(s))}{256}\left<\Delta,\Delta\right>\,.
\end{equation}

\section{Result for $L(s-9,\Delta)\,L(s-10,\Delta)$}

Combining together (\ref{eq:L(s,Delta)L(s-1,Delta)}) and (\ref{eq:part1}) we obtain the expression for the product $L(s-9,\Delta)\,L(s-10,\Delta)$:
\begin{equation}
\label{eq:L(s-9,Delta)L(s-10,Delta)}
L(s-9,\Delta)\,L(s-10,\Delta)=\frac{3\cdot 2^{13}\pi^{11}(232\,A_1(s-9)-A_2(s-9))}{(1+3\cdot2^{13-s}+2^{31-2s})\,\Gamma(s-9)}\left<\Delta,\Delta\right>\,.
\end{equation}

Now we evaluate this result in the form $L(s-9,\Delta)\,L(s-10,\Delta)=R_{\Delta}(s)\,P_{\Delta}(s)\,\left<\Delta,\Delta\right>$ for each $s\in\{12,\dotsc,19\}$ computing the rational coefficient $R_{\Delta}$ and the corresponding power of $\pi$, see Table \ref{tbl:L(s-9,Delta)L(s-10,Delta)}.  The numerical value of the Petersson inner product 
\begin{equation}
\left<\Delta,\Delta\right>=0.000001035362056
\end{equation}
is computed in the section \ref{sec:Petersson}.

\begin{table}[!tb]
\caption{$L(s-9,\Delta)\,L(s-10,\Delta)$}\label{tbl:L(s-9,Delta)L(s-10,Delta)}
$$
\renewcommand{\arraystretch}{0}
\begin{array}{|c|c@{=}c|c|c|}\hline
\rule{0pt}{11pt}s&\multicolumn{2}{c|}{R_{\Delta}}&P_{\Delta}&L(s-9,\Delta)\,L(s-10,\Delta)\\
&\multicolumn{2}{c|}{\rule{0pt}{2pt}}&&\\\hline&\multicolumn{2}{c|}{\rule{0pt}{2pt}}&&\\
\strut 12&\dfrac{32768}{225}      &\dfrac{2^{15}}{3^2\cdot5^2        }&\pi^5   &0.046143339818118\\
&\multicolumn{2}{c|}{\rule{0pt}{2pt}}&&\\\hline&\multicolumn{2}{c|}{\rule{0pt}{2pt}}&&\\
\strut 13&\dfrac{4096}{81}        &\dfrac{2^{12}}{3^4                }&\pi^7   &0.158130732552033\\
&\multicolumn{2}{c|}{\rule{0pt}{2pt}}&&\\\hline&\multicolumn{2}{c|}{\rule{0pt}{2pt}}&&\\
\strut 14&\dfrac{2048}{189}       &\dfrac{2^{11}}{3^3        \cdot7  }&\pi^9   &0.334433094416363\\
&\multicolumn{2}{c|}{\rule{0pt}{2pt}}&&\\\hline&\multicolumn{2}{c|}{\rule{0pt}{2pt}}&&\\
\strut 15&\dfrac{8192}{4725}      &\dfrac{2^{13}}{3^3\cdot5^2\cdot7  }&\pi^{11}&0.528115574483468\\
&\multicolumn{2}{c|}{\rule{0pt}{2pt}}&&\\\hline&\multicolumn{2}{c|}{\rule{0pt}{2pt}}&&\\
\strut 16&\dfrac{16384}{70875}    &\dfrac{2^{14}}{3^4\cdot5^3\cdot7  }&\pi^{13}&0.694972239760782\\
&\multicolumn{2}{c|}{\rule{0pt}{2pt}}&&\\\hline&\multicolumn{2}{c|}{\rule{0pt}{2pt}}&&\\
\strut 17&\dfrac{8192}{297675}    &\dfrac{2^{13}}{3^5\cdot5^2\cdot7^2}&\pi^{15}&0.816559651925946\\
&\multicolumn{2}{c|}{\rule{0pt}{2pt}}&&\\\hline&\multicolumn{2}{c|}{\rule{0pt}{2pt}}&&\\
\strut 18&\dfrac{8192}{2679075}   &\dfrac{2^{13}}{3^7\cdot5^2\cdot7^2}&\pi^{17}&0.895457859377812\\
&\multicolumn{2}{c|}{\rule{0pt}{2pt}}&&\\\hline&\multicolumn{2}{c|}{\rule{0pt}{2pt}}&&\\
\strut 19&\dfrac{65536}{200930625}&\dfrac{2^{16}}{3^8\cdot5^4\cdot7^2}&\pi^{19}&0.942700248523234\\
&\multicolumn{2}{c|}{\rule{0pt}{2pt}}&&\\\hline
\end{array}
$$
\end{table}

\newpage
\section{Computation of $L(s,\Delta\otimes g_{20})$}

We apply once again \cite[(2.4)]{shimura:1976} similarly as in section \ref{sec:L(s,Delta_x_G(2,2))}.  The main difference is that the Petersson inner product is taken for both modular forms being cusp forms and for the full modular group $\mathrm{SL}_2(\mathbb{Z})$:
\begin{equation}
\begin{split}
L(s,\Delta\otimes g_{20})&=\frac{(4\pi)^s}{2\,\Gamma(s)}\int_{\Phi_1}\overline{g_{20}}\,\Delta\,E_{8,1}(z,s-19)\,y^{s-1}\,dx\,dy\\
&=\frac{(4\pi)^s}{2\,\Gamma(s)}\int_{\Phi_1}\overline{g_{20}}\,\Delta\,E_{8,1}(z,s-19)\,y^{18}\,y^{s-19}\,dx\,dy\\
&=\frac{(4\pi)^s}{2\,\Gamma(s)}\left<g_{20},\Delta\,y^{s-19}\,E_{8,1}(z,s-19)\right>\\
&=\frac{(4\pi)^s}{2\,\Gamma(s)}\left<g_{20},\mathcal{H}ol\left(\Delta\,y^{s-19}\,E_{8,1}(z,s-19)\right)\right>\\
&=\frac{(4\pi)^{19}}{2\,\Gamma(s)}\left<g_{20},\mathcal{H}ol\left(\Delta\,(4\pi y)^{s-19}\,E_{8,1}(z,s-19)\right)\right>\,,
\end{split}
\end{equation}
In this case the critical values of $s$ are $12,\dotsc,19$.  We verify that $s-19\leqslant0$, $s-19+8=s-11>0$, then the series $E_{8,1}(z,s-19)$ is a nearly holomorphic modular form for all these $s\in\{12,\dotsc,19\}$.  We write the Fourier expansion of $E_{8,1}(z,s-19)$ using \cite[Proposition 2.2]{panchishkin:2003}:
\begin{equation}
\begin{split}
(4&\pi y)^{s-19}\,E_{8,1}(z,s-19)\\
&=(4\pi y)^{s-19}\,\Big[2\zeta(2s-30)\\
&\quad+\frac{(-2\pi i)^{2s-30}\,(-1)^{s-19}\Gamma(2s-31)}{(4\pi y)^{2s-31}\,\Gamma(s-11)\,\Gamma(s-19)}\,2\zeta(2s-31)\\
&\quad+\frac{2(-2\pi i)^{2s-30}(-1)^{s-19}}{\Gamma(s-11)}\times\\
&\qquad\times\sum_{n=1}^\infty\sum_{d|n}d^{2s-31}\,W(4\pi ny,s-11,s-19)\,q^n\Big]\\
%
%&=2\zeta(2s-30)(4\pi y)^{s-19} + 2(2\pi)^{2s-30}\frac{\zeta(2s-31)\,\Gamma(2s-31)}{\Gamma(s-11)\,\Gamma(s-19)}\,(4\pi y)^{12-s}\\
%&\quad+(4\pi y)^{s-19}\frac{2(2\pi)^{2s-30}}{\Gamma(s-11)}\sum_{n=1}^\infty\sum_{d|n}d^{2s-31}\,W(4\pi ny,s-11,s-19)\,q^n\\
%
&=D_0^\prime\,(4\pi y)^{12-s}+D_0^{\prime\prime}\,(4\pi y)^{s-19}\\
&\quad+\sum_{n=1}^\infty 2\,(2\pi)^{2s-30}\,\sum_{d|n}d^{2s-31}\,\frac{W(4\pi ny,s-11,s-19)}{\Gamma(s-11)}(4\pi y)^{s-11}\,q^n\,,
\end{split}
\end{equation}
where
\begin{equation}
\label{eq:D0}
\begin{split}
D_0^\prime&=D_0^\prime(s)=2\,(2\pi)^{2s-30}\frac{\Gamma(2s-31)\,\zeta(2s-31)}{\Gamma(s-11)\,\Gamma(s-19)}\,,\\
D_0^{\prime\prime}&=D_0^{\prime\prime}(s)=2\,\zeta(2s-30)\,.
\end{split}
\end{equation}

Since the result of the holomorphic projection in this case belongs to the one-dimensional space spanned by $g_{20}$, we need to compute just the first Fourier coefficient $B_1(s)$ of $\mathcal{H}ol(\cdot)=\sum_{n=1}^\infty B_n(s)q^n$ in order to express it as a multiple of $g_{20}$.  We compute it as the integral given by the Holomorphic Projection Lemma:
\begin{equation}
\begin{split}
&B_1(s)=\frac{1}{18!}\int_0^{+\infty} D_0^{\prime\prime}\,(4\pi y)^{s-19}\,\mathrm{e}^{-4\pi y}\,(4\pi y)^{18}\,d(4\pi y)\\
&\quad+\frac{1}{18!}\int_0^{+\infty} D_0^\prime\,(4\pi y)^{12-s}\,\mathrm{e}^{-4\pi y}\,(4\pi y)^{18}\,d(4\pi y)\\
&=\frac{\Gamma(s)}{18!}\,D_0^{\prime\prime} + \frac{\Gamma(31-s)}{18!}\,D_0^{\prime}
\end{split}
\end{equation}
The final expression is as following:
\begin{equation}
\label{eq:part2}
\begin{split}
L&(s,\Delta\otimes g_{20})=B_1(s)\,\frac{(4\pi)^{19}}{2\,\Gamma(s)}\left<g_{20},g_{20}\right>\\
&=\left(\frac{\Gamma(s)}{18!}\,D_0^{\prime\prime} + \frac{\Gamma(31-s)}{18!}\,D_0^{\prime}\right)\,\frac{(4\pi)^{19}}{2\,\Gamma(s)}\left<g_{20},g_{20}\right>\\
&=\frac{(4\pi)^{19}}{2\cdot 18!}\,\left(D_0^{\prime\prime} + \frac{\Gamma(31-s)}{\Gamma(s)}\,D_0^{\prime}\right)\,\left<g_{20},g_{20}\right>\,.
\end{split}
\end{equation}

Now we evaluate this result for each $s\in\{12,\dotsc,19\}$ in the form $L(s,\Delta\otimes g_{20})=R_{g_{20}}(s)\,P_{g_{20}}(s)\,\left<g_{20},g_{20}\right>$ computing the rational coefficient $R_{g_{20}}$ and the corresponding power of $\pi$, see Table \ref{tbl:L(s,Delta-times-g20)}.  The numerical value of the Petersson inner product 
\begin{equation}
\left<g_{20},g_{20}\right>=0.00000826554153165970
\end{equation}
is computed in the section \ref{sec:Petersson}.

\begin{table}[!t]
\caption{$L(s,\Delta\otimes g_{20})$}\label{tbl:L(s,Delta-times-g20)}
$$
\renewcommand{\arraystretch}{0}
\begin{array}{|c|c@{=}c|c|c|}\hline
\rule{0pt}{11pt}s&\multicolumn{2}{c|}{R_{g_{20}}}&P_{g_{20}}&L(s,\Delta\otimes g_{20})\\
&\multicolumn{2}{c|}{\rule{0pt}{2pt}}&&\\\hline&\multicolumn{2}{c|}{\rule{0pt}{2pt}}&&\\
\strut 12&\dfrac{524288}{2338875}         &\dfrac{2^{19}}{3^5   \cdot5^3\cdot7  \cdot11              }&\pi^{13}&5.380003562880315\\
&\multicolumn{2}{c|}{\rule{0pt}{2pt}}&&\\\hline&\multicolumn{2}{c|}{\rule{0pt}{2pt}}&&\\
\strut 13&\dfrac{2097152}{88409475}       &\dfrac{2^{21}}{3^8   \cdot5^2\cdot7^2\cdot11              }&\pi^{15}&5.618889612918517\\
&\multicolumn{2}{c|}{\rule{0pt}{2pt}}&&\\\hline&\multicolumn{2}{c|}{\rule{0pt}{2pt}}&&\\
\strut 14&\dfrac{4194304}{2791213425}     &\dfrac{2^{22}}{3^8   \cdot5^2\cdot7  \cdot11\cdot13\cdot17}&\pi^{17}&3.513063561721911\\
&\multicolumn{2}{c|}{\rule{0pt}{2pt}}&&\\\hline&\multicolumn{2}{c|}{\rule{0pt}{2pt}}&&\\
\strut 15&\dfrac{8388608}{97692469875}    &\dfrac{2^{23}}{3^8   \cdot5^3\cdot7^2\cdot11\cdot13\cdot17}&\pi^{19}&1.981288433718698\\
&\multicolumn{2}{c|}{\rule{0pt}{2pt}}&&\\\hline&\multicolumn{2}{c|}{\rule{0pt}{2pt}}&&\\
\strut 16&\dfrac{8388608}{1465387048125}  &\dfrac{2^{23}}{3^9   \cdot5^4\cdot7^2\cdot11\cdot13\cdot17}&\pi^{21}&1.303635536350500\\
&\multicolumn{2}{c|}{\rule{0pt}{2pt}}&&\\\hline&\multicolumn{2}{c|}{\rule{0pt}{2pt}}&&\\
\strut 17&\dfrac{2097152}{4396161144375}  &\dfrac{2^{21}}{3^{10}\cdot5^4\cdot7^2\cdot11\cdot13\cdot17}&\pi^{23}&1.072197252248449\\
&\multicolumn{2}{c|}{\rule{0pt}{2pt}}&&\\\hline&\multicolumn{2}{c|}{\rule{0pt}{2pt}}&&\\
\strut 18&\dfrac{4194304}{92319384031875} &\dfrac{2^{22}}{3^{11}\cdot5^4\cdot7^3\cdot11\cdot13\cdot17}&\pi^{25}&1.007825020916877\\
&\multicolumn{2}{c|}{\rule{0pt}{2pt}}&&\\\hline&\multicolumn{2}{c|}{\rule{0pt}{2pt}}&&\\
\strut 19&\dfrac{2097152}{461596920159375}&\dfrac{2^{21}}{3^{11}\cdot5^5\cdot7^3\cdot11\cdot13\cdot17}&\pi^{27}&0.994683426196918\\
&\multicolumn{2}{c|}{\rule{0pt}{2pt}}&&\\\hline
\end{array}
$$
\end{table}

\section{The main identity}
\label{sec:main_result}

Combining (\ref{eq:L(s-9,Delta)L(s-10,Delta)}) and (\ref{eq:part2}) into the original expression (\ref{eq:L(s,Sp(F_{12}))}) we get the following expression:
\begin{equation}
\label{eq:L(s,Sp(F_{12}))-final}
\begin{split}
L&(s,F_{12},spin)=L(s-9,\Delta)\,L(s-10,\Delta)\,L(s,\Delta\otimes g_{20})\\
&=\frac{3\cdot 2^{13}\pi^{11}}{\Gamma(s-9)}\frac{(232\,A_1(s-9)-A_2(s-9))}{(1+3\cdot2^{13-s}+2^{31-2s})}\left<\Delta,\Delta\right>\,\times\\
&\quad\times\frac{(4\pi)^{19}}{2\cdot 18!}\,\left(D_0^{\prime\prime}(s) + \frac{\Gamma(31-s)}{\Gamma(s)}\,D_0^{\prime}(s)\right)\,\left<g_{20},g_{20}\right>=\\
&=\frac{3\cdot 2^{50}\,\pi^{30}\,(232\,A_1(s-9)-A_2(s-9))}{18!\,\Gamma(s-9)\,(1+3\cdot2^{13-s}+2^{31-2s})}\,\times\\
&\quad\times\,\left(D_0^{\prime\prime}(s) + \frac{\Gamma(31-s)}{\Gamma(s)}\,D_0^{\prime}(s)\right)\,\left<\Delta,\Delta\right>\,\left<g_{20},g_{20}\right>\,,
\end{split}
\end{equation}
where $A_1(s)$ is given by (\ref{eq:A1}), $A_2(s)$ is given by (\ref{eq:A2}), $D_0^\prime$ and $D_0^{\prime\prime}$ are given by (\ref{eq:D0}).  For each critical value $s\in\{12,\dotsc,19\}$ we evaluate this expression in the form $L(s,F_{12},spin)=R(s)\,P(s)\,\left<\Delta,\Delta\right>\,\left<g_{20},g_{20}\right>$ computing the rational coefficient $R$ and the corresponding power of $\pi$, see Table \ref{tbl:L(s,Sp(F12))}.

\begin{table}[!ht]
\caption{$L(s,F_{12},spin)$}\label{tbl:L(s,Sp(F12))}
$$
\renewcommand{\arraystretch}{0}
\begin{array}{|c|c@{=}c|c|c|}\hline
\rule{0pt}{11pt}s&\multicolumn{2}{c|}{R}&P&L(s,F_{12},spin)\\
&\multicolumn{2}{c|}{\rule{0pt}{2pt}}&&\\\hline&\multicolumn{2}{c|}{\rule{0pt}{2pt}}&&\\
\strut 12&\dfrac{17179869184}{526246875}               &\dfrac{2^{34}}{3^7   \cdot5^5\cdot7  \cdot11              }&\pi^{18}&0.248251332624670\\
&\multicolumn{2}{c|}{\rule{0pt}{2pt}}&&\\\hline&\multicolumn{2}{c|}{\rule{0pt}{2pt}}&&\\
\strut 13&\dfrac{8589934592}{7161167475}               &\dfrac{2^{33}}{3^{12}\cdot5^2\cdot7^2\cdot11              }&\pi^{22}&0.888519130619814\\
&\multicolumn{2}{c|}{\rule{0pt}{2pt}}&&\\\hline&\multicolumn{2}{c|}{\rule{0pt}{2pt}}&&\\
\strut 14&\dfrac{8589934592}{527539337325}             &\dfrac{2^{33}}{3^{11}\cdot5^2\cdot7^2\cdot11\cdot13\cdot17}&\pi^{26}&1.174884717828030\\
&\multicolumn{2}{c|}{\rule{0pt}{2pt}}&&\\\hline&\multicolumn{2}{c|}{\rule{0pt}{2pt}}&&\\
\strut 15&\dfrac{68719476736}{461596920159375}         &\dfrac{2^{36}}{3^{11}\cdot5^3\cdot7^3\cdot11\cdot13\cdot17}&\pi^{30}&1.046349279390801\\
&\multicolumn{2}{c|}{\rule{0pt}{2pt}}&&\\\hline&\multicolumn{2}{c|}{\rule{0pt}{2pt}}&&\\
\strut 16&\dfrac{137438953472}{103859307035859375}     &\dfrac{2^{37}}{3^{13}\cdot5^7\cdot7^3\cdot11\cdot13\cdot17}&\pi^{34}&0.905990508529256\\
&\multicolumn{2}{c|}{\rule{0pt}{2pt}}&&\\\hline&\multicolumn{2}{c|}{\rule{0pt}{2pt}}&&\\
\strut 17&\dfrac{17179869184}{1308627268651828125}     &\dfrac{2^{24}}{3^{15}\cdot5^6\cdot7^4\cdot11\cdot13\cdot17}&\pi^{38}&0.875513015091950\\
&\multicolumn{2}{c|}{\rule{0pt}{2pt}}&&\\\hline&\multicolumn{2}{c|}{\rule{0pt}{2pt}}&&\\
\strut 18&\dfrac{34359738368}{247330553775195515625}   &\dfrac{2^{35}}{3^{18}\cdot5^6\cdot7^5\cdot11\cdot13\cdot17}&\pi^{42}&0.902464835857626\\
&\multicolumn{2}{c|}{\rule{0pt}{2pt}}&&\\\hline&\multicolumn{2}{c|}{\rule{0pt}{2pt}}&&\\
\strut 19&\dfrac{137438953472}{92748957665698318359375}&\dfrac{2^{37}}{3^{19}\cdot5^9\cdot7^5\cdot11\cdot13\cdot17}&\pi^{46}&0.937688313077777\\
&\multicolumn{2}{c|}{\rule{0pt}{2pt}}&&\\\hline
\end{array}
$$
\end{table}

\section{Numerical computation of Petersson inner product}
\label{sec:Petersson}

To compute numerically the Petersson inner product of $\Delta$ by itself and $g_{20}$ by itself we use the classical result by Rankin \cite[Theorem 5]{rankin:1952}:
\begin{equation}
\left<f_k,f_k\right> = \frac{(4\pi)^{1-k}\,(k-2)!}{\zeta(l)}\,\frac{\alpha_r}{\alpha_l+\alpha_r-\alpha_k}\,L(k-1,f_k)\,L(l,f_k)\,,
\end{equation}
where $f_k$ is the cusp form of weight $k=\{12,16,18,20,22,26\}$ of the form
\begin{equation}
f_k(z)=E_{k-12}(z)\,\Delta(z)\,
\end{equation}
and
\begin{equation}
\begin{split}
&4\leqslant r\leqslant k/2-2\,,\\
&l=k-r\,;
\end{split}
\end{equation}
$E_k$ denotes Eisenstein series
\begin{equation}
\begin{split}
&E_k(z)=\sum_{n=0}^\infty \alpha_k(n)\,q^n\,,\\
&\alpha_k(0)=1\,,\\
&\alpha_k=\alpha_k(1)=-\frac{2k}{B_k}\,,
\end{split}
\end{equation}
$B_k$ is a Bernoulli number.

For $f_k=\Delta$ we are able to use only one choice of critical value $l=8$.  To compute the numerical values $L(11,\Delta)$ and $L(8,\Delta)$ we used Dokchitser's $L$-functions Calculator \cite{ComputeL}.  In order to achieve the default precision (53 machine bits, which satisfies the functional equation to $1\text{\texttt{E-}}21$), it is necessary to input 12 Fourier coefficients in this case.  The obtained value is 
\begin{equation}
\left<\Delta,\Delta\right>=0.000001035362056804320948209596804\,,
\end{equation}
which coincides with the value given by Zagier in \cite[page 116]{zagier:1977} up to 11 digit (his method involves the direct summation of 250 first terms in $L$-series).

We used again Rankin's theorem to compute the Petersson inner product of $g_{20}$ by itself.  For the modular form of weight 20 there are three choices $l=12,14,16$.  For each choice of $l$ we computed the special value of $L(l,g_{20})$ using Dokchitser's $L$-functions Calculator.  It required to input $14$ Fourier coefficients of $g_{20}$ in order to achieve the default precision.  The obtained values are 
\begin{equation}
\begin{split}
\left<g_{20},g_{20}\right>&=0.000008265541531659702744699575969\text{~for $l=12$,}\\
\left<g_{20},g_{20}\right>&=0.000008265541531659703390644766954\text{~for $l=14$,}\\
\left<g_{20},g_{20}\right>&=0.000008265541531659703069998511729\text{~for $l=16$.}
\end{split}
\end{equation}

\section{Numerical verification}

The obtained values of $L(s,F_{12},spin)$ in section \ref{sec:main_result} can be numerically verified by using Dokchitser's $L$-functions Calculator by computing each term of the product in the righthand side of the identity (\ref{eq:L(s,Sp(F_{12}))}).  The computation of $L(s,\Delta)$ was already mentioned in the previous section.

To compute the values of $L(s,\Delta\otimes g_{20})$ for \mbox{$s\in\{12,\dotsc,19\}$} we have to determine the coefficients of this $L$-series first.  Using the identity (\ref{eq:Rankins_identity}) for $f=\Delta=\sum\tau(n)n^{-s}$ and \mbox{$g=g_{20}=\sum b(n)n^{-s}\in S_{20}$} we get
\begin{equation}
\begin{split}
L(s,\Delta\otimes g_{20})&=\sum_{n=1}^\infty A(n)n^{-s}\\
&=\sum_{d=1}^\infty d^{30-2s}\sum_{d_1=1}^\infty\tau(d_1)b(d_1)d_1^{-s}\\
&=\sum_{d,d_1\geqslant 1}d^{30}\tau(d_1)b(d_1)\,(d^2d_1)^{-s}\\
&=\sum_{n=1}^\infty\,\sum_{d:d^2|n}d^{30}\tau(\frac{n}{d^2})b(\frac{n}{d^2})\,n^{-s}\,.
\end{split}
\end{equation}
Therefore we obtain 
\begin{equation}
A(n)=\sum_{d:d^2|n}d^{30}\,\tau(\frac{n}{d^2})\,b(\frac{n}{d^2})\,.
\end{equation}
The ComputeL program requires some functional equation parameters such as $\Gamma$-factors and the weight.  These parameters can be deduced from the Hodge structures (see \cite{scholl:1990}) of $\Delta$ and $g_{20}$, namely 
\begin{equation}
 \begin{split}
  \Delta &\longrightarrow (0,11)+(11,0)\,,\\
  g_{20} &\longrightarrow (0,19)+(19,0)\,.
 \end{split}
\end{equation}
Therefore the Hodge structure of their tensor product (see \cite{yoshida:2001}) is
\begin{equation}
 \Delta\otimes g_{20} \longrightarrow (0,30)+(11,19)+(19,11)+(30,0)\,.
\end{equation}
The Deligne's rule gives \cite[page 329]{deligne:1979} in our case two $\Gamma$-factors: $\Gamma_\mathbb{C}(s)$ and $\Gamma_\mathbb{C}(s-11)$.  One can use the Gauss Duplication formula, which gives in our case four factors $\Gamma_{\mathbb{R}}(s)\,\Gamma_{\mathbb{R}}(s+1)\,\Gamma_{\mathbb{R}}(s-11)\,\Gamma_{\mathbb{R}}(s-10)$.  The weight that appears in the functional equation is the weight of the tensor product of two motives plus 1, which gives us $11+19+1=31$.  We need about 150 coefficients of $L$-series to obtain the default precision (the functional equation is satisfied with $1\text{\texttt{E-}}27$ precision).  Coefficient of $\Delta$ and $g_{20}$ are readily available in SAGE.  First few coefficients $A(n)$ are given in Table \ref{tbl:A(n)}.

\begin{table}
\caption{Fourier coefficients of $L(s,\Delta\otimes g_{20})$}\label{tbl:A(n)}
$$
\renewcommand{\arraystretch}{0}
\begin{array}{|r|r|r|r|}\hline
\rule{0pt}{10pt}n&\tau&g_{20}&A\\
\rule{0pt}{2pt}&&&\\\hline\rule{0pt}{2pt}&&&\\
\strut  1&         1&               1&                      1\\
\rule{0pt}{2pt}&&&\\\hline\rule{0pt}{2pt}&&&\\
\strut  2&       -24&             456&                 -10944\\
\rule{0pt}{2pt}&&&\\\hline\rule{0pt}{2pt}&&&\\
\strut  3&       252&           50652&               12764304\\
\rule{0pt}{2pt}&&&\\\hline\rule{0pt}{2pt}&&&\\
\strut  4&     -1472&         -316352&             1539411968\\
\rule{0pt}{2pt}&&&\\\hline\rule{0pt}{2pt}&&&\\
\strut  5&      4830&        -2377410&           -11482890300\\
\rule{0pt}{2pt}&&&\\\hline\rule{0pt}{2pt}&&&\\
\strut  6&     -6048&        23097312&          -139692542976\\
\rule{0pt}{2pt}&&&\\\hline\rule{0pt}{2pt}&&&\\
\strut  7&    -16744&       -16917544&           283267356736\\
\rule{0pt}{2pt}&&&\\\hline\rule{0pt}{2pt}&&&\\
\strut  8&     84480&      -383331840&        -44134904365056\\
\rule{0pt}{2pt}&&&\\\hline\rule{0pt}{2pt}&&&\\
\strut  9&   -113643&      1403363637&         46408678295058\\
\rule{0pt}{2pt}&&&\\\hline\rule{0pt}{2pt}&&&\\
\strut 10&   -115920&     -1084098960&        125668751443200\\
\rule{0pt}{2pt}&&&\\\hline\rule{0pt}{2pt}&&&\\
\strut 11&    534612&       -16212108&         -8667187482096\\
\rule{0pt}{2pt}&&&\\\hline\rule{0pt}{2pt}&&&\\
\strut 12&   -370944&    -16023861504&      19649522340790272\\
\rule{0pt}{2pt}&&&\\\hline\rule{0pt}{2pt}&&&\\
\strut 13&   -577738&     50421615062&     -29130483042689756\\
\rule{0pt}{2pt}&&&\\\hline\rule{0pt}{2pt}&&&\\
\strut 14&    401856&     -7714400064&      -3100077952118784\\
\rule{0pt}{2pt}&&&\\\hline\rule{0pt}{2pt}&&&\\
\strut 15&   1217160&   -120420571320&    -146571102587851200\\
% \rule{0pt}{2pt}&&&\\\hline\rule{0pt}{2pt}&&&\\
% \strut 16&    987136&     -8939761664&    1644106253837795328\\
% \rule{0pt}{2pt}&&&\\\hline\rule{0pt}{2pt}&&&\\
% \strut 17&  -6905934&    225070099506&   -1554319252561868604\\
% \rule{0pt}{2pt}&&&\\\hline\rule{0pt}{2pt}&&&\\
% \strut 18&   2727432&    639933818472&    -507896575261114752\\
% \rule{0pt}{2pt}&&&\\\hline\rule{0pt}{2pt}&&&\\
% \strut 19&  10661420&  -1710278572660&  -18233998180128777200\\
% \rule{0pt}{2pt}&&&\\\hline\rule{0pt}{2pt}&&&\\
% \strut 20&  -7109760&    752098408320&  -17676898755051110400\\
\rule{0pt}{2pt}&&&\\\hline
\end{array}
$$
\end{table}

Finally, we are able to compare the result in Table \ref{tbl:L(s,Sp(F12))} and both of its parts in Tables \ref{tbl:L(s-9,Delta)L(s-10,Delta)}, \ref{tbl:L(s,Delta-times-g20)} with the direct numerical computation.  These values and the absolute values of the difference with theoretical rational computation results are presented in Table \ref{tbl:numerical_parts} and Table \ref{tbl:numerical_final}.

\begin{table}
\caption{Numerical computation and comparison}\label{tbl:numerical_parts}
$$
\renewcommand{\arraystretch}{0}
\begin{array}{|c||c|c||c|c|}\hline
\rule{0pt}{16pt}s&{L(s-9,\Delta)\,L(s-10,\Delta)}&\displaystyle{\text{variation}\atop\text{from Table \ref{tbl:L(s-9,Delta)L(s-10,Delta)}}}&L(s,\Delta\otimes g_{20})&\displaystyle{\text{variation}\atop\text{from Table \ref{tbl:L(s,Delta-times-g20)}}}\\
\rule{0pt}{2pt}&&&&\\\hline\rule{0pt}{2pt}&&&&\\
\strut 12& 0.046143339853964& 3.58\text{\texttt{E-}}11& 5.38000356288032& 4.95\text{\texttt{E-}}15\\
\strut 13& 0.158130732674877& 1.23\text{\texttt{E-}}10& 5.61888961291852& 2.39\text{\texttt{E-}}15\\
\strut 14& 0.334433094676168& 2.60\text{\texttt{E-}}10& 3.51306356172191& 1.44\text{\texttt{E-}}15\\
\strut 15& 0.528115574893734& 4.10\text{\texttt{E-}}10& 1.98128843371870& 1.36\text{\texttt{E-}}15\\
\strut 16& 0.694972240300672& 5.40\text{\texttt{E-}}10& 1.30363553635050& 7.99\text{\texttt{E-}}16\\
\strut 17& 0.816559652560290& 6.34\text{\texttt{E-}}10& 1.07219725224845& 7.40\text{\texttt{E-}}16\\
\strut 18& 0.895457860073449& 6.96\text{\texttt{E-}}10& 1.00782502091688& 2.75\text{\texttt{E-}}15\\
\strut 19& 0.942700249255570& 7.32\text{\texttt{E-}}10& 0.99468342619692& 4.22\text{\texttt{E-}}16\\
\rule{0pt}{2pt}&&&&\\\hline
\end{array}
$$
\end{table}

\begin{table}
\caption{Numerical computation and comparison (final)}\label{tbl:numerical_final}
$$
\renewcommand{\arraystretch}{0}
\begin{array}{|c|c|c|}\hline
\rule{0pt}{16pt}s&L(s,F_{12},spin)&\displaystyle{\text{variation}\atop\text{from Table \ref{tbl:L(s,Sp(F12))}}}\\
\rule{0pt}{2pt}&&\\\hline\rule{0pt}{2pt}&&\\
\strut 12& 0.24825133281752& 1.98\text{\texttt{E-}}10\\
\strut 13& 0.88851913131006& 6.90\text{\texttt{E-}}10\\
\strut 14& 1.17488471874074& 9.13\text{\texttt{E-}}10\\
\strut 15& 1.04634928020366& 8.13\text{\texttt{E-}}10\\
\strut 16& 0.90599050923308& 7.04\text{\texttt{E-}}10\\
\strut 17& 0.87551301577209& 6.80\text{\texttt{E-}}10\\
\strut 18& 0.90246483655871& 7.01\text{\texttt{E-}}10\\
\strut 19& 0.93768831380622& 7.28\text{\texttt{E-}}10\\
\rule{0pt}{2pt}&&\\\hline
\end{array}
$$
\end{table}

%\newpage
\bibliographystyle{alpha}
\bibliography{L-Spin-F12}

\end{document}